\documentclass[reqno,11pt]{amsart}
\usepackage{graphics, color, amsmath, amsfonts, amssymb, amsthm, amscd}
\usepackage[colorlinks=true, pdfstartview=FitV, linkcolor=blue,
  citecolor=blue, urlcolor=blue,pagebackref=false]{hyperref}
\usepackage{mathabx}

\def\writefig#1 #2 #3 {\rlap{\kern 62.5 truemm\kern #1 truecm
    \Raise #2 truecm \hbox{#3}}}

\numberwithin{equation}{section}

\DeclareMathSymbol{\leqslant}{\mathalpha}{AMSa}{"36} 
\DeclareMathSymbol{\geqslant}{\mathalpha}{AMSa}{"3E} 
\DeclareMathSymbol{\eset}{\mathalpha}{AMSb}{"3F}     
\renewcommand{\leq}{\;\leqslant\;}                   
\renewcommand{\geq}{\;\geqslant\;}                   

\DeclareMathOperator*{\argmax}{arg\,max}       

\newcommand{\argmaxone}[1]{\argmax_{\substack{#1}}} 
\newcommand{\maxone}[1]{\max_{\substack{#1}}} 
\newcommand{\supone}[1]{\sup_{\substack{#1}}} 
\newcommand{\infone}[1]{\inf_{\substack{#1}}} 

\makeatletter
\def\captionfont@{\footnotesize}
\def\captionheadfont@{\scshape}
\long\def\@makecaption#1#2{%
  \vspace{2mm}
  \setbox\@tempboxa\vbox{\color@setgroup
    \advance\hsize-6pc\noindent
    \captionfont@\captionheadfont@#1\@xp\@ifnotempty\@xp
        {\@cdr#2\@nil}{.\captionfont@\upshape\enspace#2}%
    \unskip\kern-6pc\par
    \global\setbox\@ne\lastbox\color@endgroup}%
  \ifhbox\@ne 
    \setbox\@ne\hbox{\unhbox\@ne\unskip\unskip\unpenalty\unkern}%
  \fi
  \ifdim\wd\@tempboxa=\z@ 
    \setbox\@ne\hbox to\columnwidth{\hss\kern-6pc\box\@ne\hss}%
  \else 
    \setbox\@ne\vbox{\unvbox\@tempboxa\parskip\z@skip
        \noindent\unhbox\@ne\advance\hsize-6pc\par}%
\fi
  \ifnum\@tempcnta<64 
    \addvspace\abovecaptionskip
    \moveright 3pc\box\@ne
  \else 
    \moveright 3pc\box\@ne
\nobreak
\vskip\belowcaptionskip
\fi
\relax
}
\makeatother

\def\writefig#1 #2 #3 {\rlap{\kern #1 truecm
\raise #2 truecm \hbox{#3}}}

\newtheorem{theorem}{Theorem}[section]

\newtheorem{corollary}[theorem]{Corollary}
\newtheorem{lem}[theorem]{Lemma}
\newtheorem{prop}[theorem]{Proposition}
\newtheorem{thm}[theorem]{Theorem}

\newcommand{\cA}{\ensuremath{\mathcal A}}
\newcommand{\cB}{\ensuremath{\mathcal B}}

\newcommand{\cD}{\ensuremath{\mathcal D}}

\newcommand{\cH}{\ensuremath{\mathcal H}}
\newcommand{\cI}{\ensuremath{\mathcal I}}

\newcommand{\cL}{\ensuremath{\mathcal L}}

\newcommand{\cP}{\ensuremath{\mathcal P}}

\newcommand{\cX}{\ensuremath{\mathcal X}}

\newcommand{\cZ}{\ensuremath{\mathcal Z}}

\newcommand{\bbE}{{\ensuremath{\mathbb E}} }

\newcommand{\bbL}{{\ensuremath{\mathbb L}} }
\newcommand{\bbM}{{\ensuremath{\mathbb M}} }

\newcommand{\bbP}{{\ensuremath{\mathbb P}} }
\newcommand{\bbQ}{{\ensuremath{\mathbb Q}} }
\newcommand{\bbR}{{\ensuremath{\mathbb R}} }

\newcommand{\bbZ}{{\ensuremath{\mathbb Z}} }

\newcommand{\ga}{\alpha}
\newcommand{\gb}{\beta}
\newcommand{\gga}{\gamma}            
\newcommand{\gd}{\delta}
\newcommand{\gep}{\epsilon}           

\newcommand{\gD}{\Delta}

\newcommand{\gO}{\Omega}
\newcommand{\gl}{\lambda}
\newcommand{\gL}{\Lambda}
\newcommand{\gs}{\sigma}
\newcommand{\gS}{\Sigma}
\newcommand{\gt}{\theta}
\newcommand{\gGa}{\Gamma}
\newcommand{\gi}{\iota}
\newcommand{\gT}{\Theta}

\newcommand{\norm}[1]{|#1|}

\newcommand{\normsup}[1]{\|#1\|_{{\scriptscriptstyle\infty}}}

\def\1{\ifmmode {1\hskip -3pt \rm{I}}
\else {\hbox {$1\hskip -3pt \rm{I}$}}\fi} 

\newcommand{\ra}{\rightarrow}
\newcommand{\Ra}{\Rightarrow}

\newcommand{\lb}{\left(}
\newcommand{\rb}{\right)}
\newcommand{\lbr}{\left\{}
\newcommand{\rbr}{\right\}}

\newcommand{\lsp}{\left [}           
\newcommand{\rsp}{\right ]}          
\newcommand{\lpr}{\left.}
\newcommand{\rpr}{\right.}
\newcommand{\rabs}{\right|}
\newcommand{\labs}{\left|}

\newcommand{\wt}{\widetilde}
\newcommand{\ol}{\overline}

\newcommand{\inv}[1]{\tfrac{1}{#1}}
\newcommand{\minwith}{\land}
\newcommand{\maxwith}{\lor}
\newcommand{\graph}{\mbox{graph}}
\newcommand{\linear}{\mbox{linear}}
\newcommand{\For}{\mbox{for}}

\newcommand{\be}{\begin{equation}}
\newcommand{\ee}{\end{equation}}
\newcommand{\pr}{\prime}
\newcommand{\bydef}{\triangleq}
\newcommand{\wh}{\widehat}
\newcommand{\wc}{\widecheck}
\newcommand{\eqd}{{\buildrel d \over =}}

\newcommand{\forcenewline}{$ $}

\begin{document}
\title[\,]{
Directed Polymers in Random Environment with Heavy Tails
}
\author[A. Auffinger]{Antonio Auffinger}
\address{A. Auffinger\\
  Courant Institute of the Mathematical Sciences\\
  New York University\\
  251 Mercer Street\\
  New York, NY 10012, USA}
\email{auffing@cims.nyu.edu}

\author[O. Louidor]{Oren Louidor}
\address{O. Louidor\\
  Courant Institute of the Mathematical Sciences\\
  New York University\\
  251 Mercer Street\\
  New York, NY 10012, USA}
\email{louidor@cims.nyu.edu}

\setcounter{page}{1}
\maketitle
\begin{abstract}
  We study the model of Directed Polymers in Random Environment in $1+1$ dimensions, where the distribution at a site has a tail which decays regularly polynomially with power $\ga$, where $\ga \in (0,2)$. After proper scaling of temperature $\gb^{-1}$, we show strong localization of the polymer to a favorable region in the environment where energy and entropy are best balanced. We prove that this region has a weak limit under linear scaling and identify the limiting distribution as an $(\ga, \gb)$-indexed family of measures on Lipschitz curves lying inside the $45^{\circ}$-rotated square with unit diagonal. In particular, this shows order $n$ transversal fluctuations of the polymer. If, and only if, $\ga$ is small enough, we find that there exists a random critical temperature below which, but not above, the effect of the environment is macroscopic.
The results carry over to $d+1$ dimensions for $d>1$ with minor modifications.
\end{abstract}

\section{Introduction}
\label{sec:Intro}
Originally introduced in \cite{HUSE-HENLEY}, Directed Polymers in Random Environment is a model for an interaction between a {\it polymer chain} and a {\it medium} with microscopic impurities. The reader can find mathematical and physical background material in surveys  \cite{COMETS-ET-AL} and \cite{FISHER-HUSE}. In this model, the  medium with defects or {\it environment} is represented by a random measure $\gs$ with support in $\bbZ_+ \times \bbZ^d$, the shape of the polymer (random due to thermal fluctuations) is represented by the $Z_+ \times \bbZ^d$-trajectory of a simple random walk on $\bbZ^d$ and the interaction is expressed as a measure change. More precisely, conditioned on $\gs$, the $n$-monomer polymer chain is a random nearest-neighbor path $s:[0,n] \cap \bbZ \to \bbZ^d$ chosen according to the following Gibbs distribution:
\begin{equation}
\label{eqn:MuInIntro}
\mu_{\gb}^{\gs}(s) =
  \frac{1}{Q_{\gb}^{\gs}} \, \exp \lb -\gb \cH^{\gs}(s) \rb
    \quad ; \qquad \cH^{\gs}(s) = -\gs(s),
\end{equation}
where $\cH^{\gs}(s)$ is the Hamiltonian or {\it energy} of $s$, $\gs(s)$ is to be understood as the measure under $\gs$ of $\graph(s) = \{(x, s(x)) :\, x \in [0,n] \cap \bbZ \}$, the parameter $\gb \in [0,\infty)$ represents the overall strength of the interaction (inverse temperature) and $Q_{\gb}^{\gs}$ is the normalizing constant.

It is believed \cite{FISHER-HUSE} that there exist thermodynamic phases in which the effects of the environment
on the shape of the polymer become macroscopic in scale. More explicitly, it is expected that in $d \leq 2$ and any finite temperature, or $d > 2$ and low enough (but not necessarily zero, at least for $d$ small enough) temperature, the behavior of the polymer path is super-diffusive with transversal fluctuations of order $n^{1/2+\gep}$ for some $\gep > 0$, while in infinite ($d \leq 2$) or high enough ($d >2$) temperature the polymer is diffusive with fluctuations of order $n^{1/2}$.
These two phases - the {\it strong-disorder} phase and 
{\it weak-disorder} phase (respectively) and the conditions for their occurrence
should be universal with respect to a large class of environment distributions, including the case where $\{\gs(\{z\}) :\: z \in \bbZ_+ \times \bbZ^d\}$ are i.i.d. with distribution whose tail decays sufficiently fast. Partial results in this direction have been established in \cite{BOLTHAUSEN, COMETS-ET-AL, SPENCER, PIZA}.

In the strong disorder phase the attraction to regions in the environment with relatively low energy has a non-negligible counter-effect to the entropic tendency of the polymer to diffuse, resulting in a {\it pinning} or {\it localization} of the polymer to a region in the environment where the balance between {\it entropy} and energy is optimal. The transversal fluctuations of the polymer are therefore significantly effected by the fluctuations in the shape of this {\it favorable} region. This pinning becomes absolute in the extreme case of zero temperature (formally, the weak limit of $\mu_{\gb}^{\gs}$ as $\gb \to \infty$). Here entropy no longer plays a part and the system is uniformly in one of its {\it ground states} - states in which the polymer path is a minimizer of the energy $\cH^{\gs}(\cdot)$ (equivalently, maximizer of $\gs(\cdot)$). Universal super-diffusive behavior is expected here as well and some progress has been made in this direction (see \cite{BAIK-DEIFT, JOHANSSON, LICEA-NEWMAN-PIZA, NEWMAN-PIZA} and also the surveys on the equivalent models of Directed First/Last Passage Percolation in \cite{KESTEN, MARTIN}).

In this paper we address the case where $(\gs(\{z\}))_z$ are i.i.d. but their  distribution has a (right) tail which is heavy enough to fall outside the universality classes discussed above. Inspired by \cite{Ha07}, we assume that the tail of the distribution at each site is regularly varying with index $\ga \in (0,2)$, namely
\begin{equation}
\label{eqn:HeavyTailCondition}
\bbP(\gs(\{z\}) > t) = t^{-\ga} L(t) \,,
\end{equation}
where $L(t)$ is a slowly varying function ($L(tu)/L(t) \to 1$ as $t \to \infty$ for all $u > 0$ - this, of course, includes all constants). We also assume that 
$\gs(\{z\})$ are positive absolutely-continuous random variables and, for the sake of simplicitly, treat only the $d=1$ case. We show [Theorem~\ref{thm:QuenchedLocalization}] that for any finite temperature, the polymer is localized to the path along which the energy is minimal, i.e. thermal fluctuations are negligible - this is because entropy is of smaller order compared to energy in this case. Consequently, it is natural to let $\gb$ go to zero with $n$ and indeed if
\begin{equation}
\label{eqn:GbCondIntro}
\gb_n = \gb n^{1-2/\alpha} L_0(n)
\end{equation}
where $L_0(n)$ is a related slowly varying function, the system exhibits a non-trivial interplay between energy and entropy. Under this scaling, we show [Theorem~\ref{thm:QuenchedLocalization}] that the $n$-monomer polymer chain is localized
(in probability, exponentially fast) to a cylindrical region of diameter $o(n)$ around a random {\it favorable curve} [Equation \eqref{Eqn.FiniteGgaHatDef}], i.e. one that optimally balances entropy and energy under $\gs$. Zero (resp. infinite) temperature behavior is recovered if $\gb_n$ grows faster (resp. slower) than $n^{1-2/\alpha} L_0(n)$.

A weak limit for the distribution of the favorable curve under linear scaling
is then shown to exist [Theorem~\ref{thm:ScalingLimit}] and the limiting distribution (on a proper topological space of directed curves lying in the $45^{\circ}$-rotated square with unit diagonal) is explicitly described.
The limit is constructed as the distribution of the (almost-surely) unique global solution to a variational problem on the space of curves [Equation \eqref{Eqn.GgaHatDef}]. The functional being maximized is random and can be viewed as assigning to a curve the difference between its entropy gain and energy cost under a proper (random) limit environment. These limiting distributions form a two-parameter family of measures on curves $\{\bbM_{\ga, \gb}$, $\ga \in (0,2), \gb \in [0,\infty]\}$ [Proposition~\ref{prop:MIsDistingushiable}], where $\bbM_{\ga, \infty}$ is the scaling limit of the ground-state path, which was studied in \cite{Ha07}.

Combining these results we obtain [Corollary~\ref{cor:ScalingLimitOfPolymerPath}] a linear scaling limit for the unconditional distribution of the polymer's path
and transversal fluctuations of order $n$ - quite different from the light tail case. This happens for all $\gb > 0$ (under scaling \eqref{eqn:GbCondIntro}). Nevertheless [Proposition~\ref{prop:PhaseTransitionAtZero}], when $\ga$ is small enough there exists a random variable $\gb_c$ positive a.s. (but arbitrarily small with positive probability) such that if $\gb < \gb_c$ then the polymer localizes around the $x$-axis ($s \equiv 0$), thereby exhibiting an infinite-temperature behavior and if $\gb > \gb_c$, the polymer tends to drift away from the axis. This can be viewed as a {\it quenched phase transition} with a random threshold value.

Although we only treat the $1+1$ dimensional case here, the problem is exactly the same for all dimensions, with minor changes due to the different geometry. In particular, for any $d$, the right normalization is $\gb_n = \gb n^{1-(d+1)/\ga} L_0(n)$ and the limit curves live inside a $d+1$ regular polyhedron. 

\section{Setup and Results}
\label{Edt.Sec.SetupAndResults}
In this section we introduce notation and state our main results.
Let $n$ be a positive natural number. Set $\bbL = \{a (1,1) + b (1,-1) :\: a, b \in \bbZ \}$ and $\bbL_n = \frac{1}{n} \bbL$.
Define $\cD = \{(x,y) \, :\; |y| \leq min(x, 1-x)\}$,
$\cD^0 = \cD \setminus \{(0,0), (1,0)\}$ and $\cD_n$, $\cD^0_n$ as the intersection of $\cD$, $\cD^0$ with $\bbL_n$. Let $\cL^0$ be
the set of all real Lipschitz functions on $[0,1]$ with
Lipschitz constant 1, vanishing at $0$ and $1$ and 
\[
    \cL^0_n = \lbr  s \in \cL^0 : \: (k/n, s(k/n)) \in \bbL_n ,\, k=0, \dots, n \rbr.
\]
$\cL^0_n$ can be viewed as
the set of linearly interpolated $1/n$-scaled trajectories
of a simple random walk conditioned to hit $0$ at time $n$ and
hence a finite subset of $\cL^0$.
Note that $\gga \in \cL^0 \; \Ra \; \mbox{graph}(\gga) \; \subseteq \cD$.
Endow $\cL^0$, $\cL^0_n$ with the $L^{\infty}$ norm and Borel sigma algebras
$\cB^0$ and $\cB^0_n$.
Let $\cP$ (resp. $\cP_n$) be the set of all probability measures
on $\cL^0, \cB^0$ (resp. $\cL^0_n, \cB^0_n$).
We shall treat $\cP_n$ as a subset of
$\cP$.

The entropy of a $\gga \in \cL^0$ curve is $-E(\gga)$ where
\begin{equation}
\label{eqn:EntropyDef}
E(\gga) = \int_0^1 e(\gga^{\prime}(x)) d x
\end{equation}
and $e:[-1,1] \to \bbR$ is defined as
\begin{equation}
\label{eqn:LittleEDef}
e(x) = \inv{2} [(1+x) \log(1+x) + (1-x) \log(1-x)]\,.
\end{equation}
$E$ is well defined since $\gga$ is differentiable almost everywhere with $|\gga^{\prime}(x)| \leq 1$. In fact, it is the rate function in the large deviations principle for the sequence of uniform measures on $\cL^0_n$ (this is essentially Mogulskii's Theorem - see ~\cite{DZ}, Section 5.1).

$\gs_n$ will denote the $\inv{n}$-scaled and $\cD^0_n$-restricted version of  $\gs$. That is, $\gs_n$ is a positive measure on $\cD^0_n$ with i.i.d weights which satisfy \eqref{eqn:HeavyTailCondition}. The scaled analog of $\mu^\gs_{\gb}$ is
\begin{equation}
\mu_{n, \gb}(s) = 
  \frac{1}{Q_{n, \gb}} \, \exp \lb \gb \gs_n(s) \rb
  \quad \For \ s \in \cL^0_n.
\end{equation}

It is a standard fact (see Section 1.1 in \cite{Resnick}) that the 
distribution of $\gs_n(\{z\})$ is in the max-domain of attraction of the Fr\'{e}chet distribution, namely there exist $(b_n)_{n \geq 1}$ such that if $(U^i_n, Z^i_n) \in (\bbR_+, \cD^0_n)$ are the value and position of the non-ascending $i$-th order statistic of $(\gs_n(\{z\}))_z$, then for any fixed $k \geq 1$
\begin{equation}
\label{eqn:KMaxsLimit}
((b_n^{-1} U^i_n,\: Z^i_n))_{i=1}^k \Ra ((V^i,\: Z^i))_{i=1}^k\,,
\end{equation}
as $n \to \infty$, where the limit is non-degenerate. The constants $b_n$ can be written as
\begin{equation}
\label{def:bnNiceForm1}
b_n = n^{2/\ga} L_0(n)
\end{equation}
where $L_0(n)$ is a related slowly varying function.
This is the motivation behind the scaling of the temperature and accordingly if $\gb_n$ denotes inverse temperature at system order $n$, we assume
\begin{equation}
\label{eqn:GbNConditions}
\lim_{n \to \infty}  \frac{b_n}{n} \gb_n= \ol{\gb}_{\infty},
\end{equation}
where $\ol{\gb}_{\infty} \in [0,\infty] = [0,\infty) \cup \{\infty\}$. This is a more explicit version of \eqref{eqn:GbCondIntro}.

The favorable curve, around which concentration occurs, is $\gga_{n, \gb_n}^*$, where 
\begin{equation}
\label{Eqn.FiniteGgaHatDef}
\gga_{n, \gb}^* = \argmaxone{\gga \in \cL^0} \lb \gb \gs_n(\gga) - n E(\gga)\rb .
\end{equation}
The properness of this definition is discussed in Lemma~\ref{lem1}. Localization is ``exponentially fast in probability'', by which we mean
\begin{thm}
\label{thm:QuenchedLocalization}
For all $\gep > 0$, $\gd > 0$
there exist $\nu > 0$ such that
\[
\begin{array}{ll}
\mu_{n, \gb_n} \lb \normsup{s-\gga_{n, \gb_n}^*} > \gd \rb \leq e^{-n \nu}
    \qquad \qquad & \text{if } \ol{\gb}_{\infty} < \infty \,,\\
\mu_{n, \gb_n} \lb \normsup{s-\gga_{n, \gb_n}^*} > \gd \rb \leq e^{-\nu b_n \gb_n }
    & \text{if } \ol{\gb}_{\infty} = \infty\,,
\end{array}
\]
with $\bbP_n$-probability at least $1-\gep$ as long as $n$ is large enough.
\end{thm}

Denote by $\bbM_{n, \gb} \in \cP$ the distribution of $\gga^*_{n, \gb}$, namely
$\bbM_{n, \gb}(\cdot) = \bbP_n(\gga^*_{n, \gb} \in \cdot)$.
Under \eqref{eqn:GbNConditions} the sequence of measures
$(\bbM_{n, \gb_n})_{n \geq 1}$ has a weak limit in $\cP$. This limit measure is constructed on top of the scaling limit of the position and weight of the environment point masses, i.e. the infinite collection $((V^i, Z^i))_{i=1}^{\infty}$ of which each finite subset $((V^i, Z^i))_{i=1}^k$ has a law as in the limit in \eqref{eqn:KMaxsLimit}. It is a standard fact that
$(Z^i)_{i=1}^{\infty}$, $(V^i)_{i=1}^{\infty}$ are independent of each other and
\begin{equation}
\label{eqn:ExtremesLimit}
Z^i \sim \text{Uniform}(\cD)  \quad \text{i.i.d.} \quad ; \qquad V^i \eqd T_i^{-\inv{\ga}}
\end{equation}
where $T_i$ is the sum of $i$ independent exponentials with rate $1$.
We can then define the ``limit environment'' as
\[
  \pi_{\infty} = \sum_i V^i \gd \lb \cdot - Z^i \rb \,.
\]
Note that while $\pi_{\infty}$ as a
measure may be infinite (for $\ga \geq 1$), as a function $\pi_{\infty}(\gga) = \pi_{\infty}(\graph(\gga))$, it is bounded on $\cL^0$,
with $\bbP_{\infty}$-probability $1$, where we denote by $\bbP_{\infty}$ the underlying measure. This follows from Theorem 2.1 in \cite{Ha07}.

The limit curve $\wh{\gga}_{\infty, \gb}$ is defined, analogously to the finite case, as the solution to a variational problem on $\cL^0$, namely
\begin{equation}
\label{Eqn.GgaHatDef}
\wh{\gga}_{\infty, \gb} =
    \lbr \begin{array}{ll}
        \argmaxone{\gga \in \cL^0} \, \lb \gb \pi_{\infty}(\gga) - E(\gga)\rb
            \quad   & \gb < \infty  \\
        \argmaxone{\gga \in \cL^0} \, \pi_{\infty}(\gga)
                    & \gb = \infty
    \end{array} \rpr
\end{equation}
A maximizer always exists and it is unique. This will be proved in Lemma~\ref{lem1}. Formally, however, we shall set $\wh{\gga}_{\infty, \gb} \equiv \infty$, if one of these conditions fails.
With $\bbM_{\ga, \gb}(\cdot) = \bbP_{\infty} (\wh{\gga}_{\infty, \gb} \in \cdot)$, we can now state
\begin{thm}
\label{thm:ScalingLimit}
$\bbM_{n, \gb_n} \Ra \bbM_{\ga, \ol{\gb}_{\infty}}$ as $n \to \infty$ in $\cP$.
\end{thm}

The following is an immediate corollary of Theorems~\ref{thm:QuenchedLocalization} and ~\ref{thm:ScalingLimit}.
Let $S_{n,\gb_n}$ be a random variable, taking values in $\cL^0_n$, such that
conditioned on the environment $\gs_n$, its distribution is
$\mu_{n,\gb_n}$.
The unconditional distribution of
$S_{n, \gb_n}$, which we will denote by $\bbQ_{n,\gb_n}$ is obtained by averaging over the environment, namely $\bbQ_{n,\gb_n} = \bbP_n \mu_{n, \gb_n}$. Then
\begin{corollary}
\label{cor:ScalingLimitOfPolymerPath}
$\bbQ_{n, \gb_n} \Ra \bbM_{\ga, \ol{\gb}_{\infty}}$ as $n \to \infty$ in $\cP$.
\end{corollary}

In order to justify that these localization results indeed imply a qualitative $\gT(n)$ change in the shape of the polymer, we have to argue, in addition, that with positive
$\bbP_{\infty}$-probability $\wh{\gga}_{\infty, \gb} \not\equiv 0$. This is included in the next proposition. For two random variables $X$, $Y$ we write $X \lneqq_s Y$ if $X \le_s Y$ but not $X \eqd Y$.
\begin{prop}
\label{prop:MIsDistingushiable}
If  $ 0 \leq \gb_2 < \gb_1 \leq \infty$ then
$E(\wh{\gga}_{\infty, \gb_2}) \lneqq_s E(\wh{\gga}_{\infty, \gb_1})$, where
$\wh{\gga}_{\infty, \gb_i}$ is distributed according to $\bbM_{\ga, \gb_i}$.
In particular $\bbM_{\ga, \gb_1} \neq \bbM_{\ga, \gb_2}$.
\end{prop}
Hence, $\bbM_{\ga, \gb}$ for $\gb \in (0, \infty)$ is different from both
$\bbM_{\ga, 0} = \gd_0$ - the Dirac-mass on the zero function $0 \in \cL^0$ and $\bbM_{\ga, \infty}$ - the distribution of the {\it last passage path}, i.e. the path along which the sum of the mass weights of $\pi_{\infty}$ (thought of as passage times) is maximal. The latter was studied in \cite{Ha07} as $P^*$.

Nevertheless, it is still quite possible that there exists a non-degenerate random $\pi_{\infty}$-dependent $\gb_c$ such that $\wh{\gga}_{\infty, \gb} \equiv 0$ if $\gb < \gb_c$, but $\wh{\gga}_{\infty, \gb} \not \equiv 0$ if $\gb > \gb_c$. This will show a (random) {\it phase-transition}-like phenomenon, where depending on whether the temperature is below or above a random threshold, the effect of the environment is microscopic or macroscopic, measured on the scale of $\gT(n)$. To make this precise, let us define
\[
\wh{w}_{\infty, \gb} = \maxone{\gga \in \cL^0} \, \lb \gb \pi_{\infty}(\gga) - E(\gga) \rb
\quad \text{and} \quad
\gb_c = \inf \{\gb \geq 0: \wh{w}_{\infty, \gb} > 0\ \}\,.
\]
The following proposition shows that this indeed occurs for $\ga$ small enough. The restrictions on $\ga$ are not sharp.
\begin{prop}
\forcenewline
\label{prop:PhaseTransitionAtZero}
\begin{enumerate}
\item
	\label{item:PhaseTransitionAtZero0}
	$\gb_c$ is well defined $\bbP_{\infty}$-a.s.
\item
\label{item:PhaseTransitionAtZero1}
    $\wh{\gga}_{\infty, \gb} \equiv 0$ if $\gb < \gb_c$ and
    $\wh{\gga}_{\infty, \gb} \not \equiv 0$ if $\gb > \gb_c$.
\item
\label{item:PhaseTransitionAtZero2}
    if $\alpha \in [\inv{2}, 2)$ then $\gb_c = 0$ with $\bbP_{\infty}$-probability $1$.
\item
\label{item:PhaseTransitionAtZero3}
    if $\alpha \in (0, \inv{3})$ then $\gb_c > 0$ with $\bbP_{\infty}$-probability $1$.
\end{enumerate}
\end{prop}

\subsection{Organization of the paper}
In the remainder of the text, we prove the results in this section. Section~\ref{Outline.Sec.Proofs}
contains some preliminary definitions and tools, on which we base our proofs. Section~\ref{sub:localization} contains proofs for Theorems~\ref{thm:QuenchedLocalization},\ref{thm:ScalingLimit} and Corollary~\ref{cor:ScalingLimitOfPolymerPath}. Section~\ref{Section:PropLimitDist} contains proofs for Propositions~\ref{prop:MIsDistingushiable} and ~\ref{prop:PhaseTransitionAtZero}. Finally, Section~\ref{sub:DeferredProofs} contains proofs for some of the results in Section~\ref{Outline.Sec.Proofs}, which we deferred.

\section{Preliminaries}
\label{Outline.Sec.Proofs}

\subsection{Environments and $\cL^0$ Curves}
\label{Outline.Sub.EnvironmentsAndLipCurves}
We shall call an {\it environment} any positive (possibly infinite) measure on $\cD^0$ with countable support, for which the collection of mass weights can be ordered in non-ascending order. If $\gs$ is an environment, we denote by $(v_{\gs}^i,\, z_{\gs}^i = (x_{\gs}^i, y_{\gs}^i)) \in \bbR_+ \times \cD^0$ the position and weight of the $i$-th mass in this order (if there are masses with equal weight, we suppose that they come with a prescribed order). To the collection $\{(v_{\gs}^i, z_{\gs}^i)\}_{i=1}^{|\gs|}$, where $|\gs|$
is the cardinality of the support of $\gs$, we add two pairs
$(v_{\gs}^0, z_{\gs}^0), (v_{\gs}^{\infty}, z_{\gs}^{\infty})$ with
$v_{\gs}^0 = v_{\gs}^{\infty} = 0$ and $z_{\gs}^0 = (0,0)$, $z_{\gs}^{\infty} = (1,0)$. This will simplify things later on. Thus, $\gs$ is identified with two sequences:
\[
v_{\gs} = (v_{\gs}^i) : \: i=0, 1, \dots , |\gs|, \infty)
	\ , \quad
z_{\gs} = (z_{\gs}^i) : \: i=0, 1, \dots , |\gs|, \infty)
\]
and
\[
\gs = \sum_{i=0,1,\dots,|\gs|, \infty} v_{\gs}^i \gd (\cdot - {z_{\gs}^i}).
\]
Both $\gs_n$ and $\pi_{\infty}$ in Section~\ref{Edt.Sec.SetupAndResults}
are environments under this definition.

Given $i \neq j \in \{0,1, \dots, |\gs|, \infty\}$, let $\gD x_{\gs}(i,j) = x_{\gs}^j - x_{\gs}^i$, $\gD y_{\gs}(i,j) = y_{\gs}^j - y_{\gs}^i$, $a_{\gs}(i,j) = \gD y_{\gs}(i,j) / \gD x_{\gs}(i,j)$ and set
\[
	\cI(\gs) = \{ \gi \subseteq \{0,1, \dots, |\gs|, \infty\} : \quad
    0, \infty \in \gi ,\,\, |a_{\gs}(i,j)| \leq 1 \ \forall i \neq j \in \gi\}.
\]
If $\gi \in \cI(\gs)$ is finite, we shall also treat it as a sequence of indices
$\gi = (\gi_j)_{1 \leq j \leq |\gi|}$ ordered according to the $x$-coordinate of the indexed point (i.e. $j < k \Ra x_{\gs}^{\gi_j} < x_{\gs}^{\gi_k}$). $m(\gs)$ will denote the ``mesh'' of $\gs$, defined as
\[
m(\gs) = \inf_{i \neq j} \lb |\gD x_{\gs}(i,j)| \minwith ||a_{\gs}(i,j)|-1| \rb,
\]
where the infimum is over all $i, j \in \{0,1, \dots, |\gs|, \infty\}$.
The distance between two environments $\gs$, $\wt{\gs}$ with equal cardinality
$|\gs| = |\wt{\gs}| = \chi \leq \infty$ is given by
\[
d(\gs, \wt{\gs}) = d_{\chi}(\gs, \wt{\gs}) =
\normsup{v_{\gs} - v_{\wt{\gs}}} \maxwith \normsup{z_{\gs} - z_{\wt{\gs}}}.
\]
This is a well-defined metric on $\gS_{\chi} -$ the space of all environments with cardinality $\chi$. We shall also use $\gS_{\chi, n} \subseteq \gS_{\chi}$ for the subset of environments supported on $\cD_n^0$ and $\gS_{\chi}(m_0)  \subseteq \gS_{\chi}$ for the subset of environments $\gs$ with $m(\gs) > m_0$. The intersection is denoted by $\gS_{\chi,n}(m_0)$.

Given $\phi: \cX \to \bbR$, where $\cX$ is a closed subset of $[0,1]$ with $L=\min (\cX)$ and $R=\max(\cX)$, we define $\linear(\phi)$ as the $[L,R] \to \bbR$ function obtained from $\phi$ by linearly interpolating inside all intervals $\{[l(x), r(x)]  :\: x \notin \cX\}$, where $l(x) = \max (\cX \cap [0,x])$ and $r(x) = \min (\cX \cap [x,1])$. Given $\cZ \subseteq \cD$, which is a graph of such function $\phi_{\cZ}$, we set 
$\linear(\cZ) = \linear(\phi_{\cZ})$. 

The following two mappings between $\cL^0$ and $\cI(\gs)$ will be used often in what follows. 
$I_{\gs}: \cL^0 \ra \cI(\gs)$ is defined as
\[
I(\gga) = \lbr j \in \{0, 1, \dots, |\gs|, \infty \} : \:
    z_{\gs}^j \in \graph(\gga) \rbr   \quad \text{for} \: \gga \in \cL^0,
\]
and $\gGa_{\gs}: \cI(\gs) \ra \cL^0$ as
\[
\gGa_{\gs}(\gi) = \linear \lb \overline{\{ z_{\gs}^j : \: j \in \gi \}} \rb \quad \text{for} \: \gi \in \cI(\gs).
\]
The validity of these definitions is not difficult to verify.
Finally, define $\gL_{\gs} : \cL^0 \to \cL^0$ as
\[
	\gL_{\gs}(\gga) = \gGa \lb I_{\gs}(\gga) \rb \quad \text{for} \: \gga \in \cL^0.
\]

Below are three technical propositions which we use later in the text. We defer their proofs to Section~\ref{sub:DeferredProofs}. Recall the definition of $E$ in \eqref{eqn:EntropyDef}.
\begin{prop}
\label{prop:EProps}
\forcenewline
\begin{enumerate}
\item
	\label{item:EProps1}
	$E$ is lower semi-continuous, strictly convex and positive away from $0$.
\item
	\label{item:EProps2}
	If $\cZ$ is a closed subset of $\graph(\gga)$ for $\gga \in \cL^0$, then
	$E(\linear(\cZ)) \leq E(\gga)$. In particular, $E(\gL_{\gs}(\gga)) \leq E(\gga)$ for any $\gga \in \cL^0$ and any environment $\gs$.
\end{enumerate}
\end{prop}
\begin{prop}
\label{Thm.Prop.Continuity}
\forcenewline
\begin{enumerate}
\item
    \label{Item.IContinuity}
    For all $\chi_0 < \infty$, $\gs \in \gS_{\chi_0}$ there exists $d_0 = d_0(\gs) > 0$ such that
    $\cI(\wt{\gs}) \subseteq \cI(\gs)$ for all $\wt{\gs} \in \gS_{\chi_0}$ with
    $d(\gs, \wt{\gs}) < d_0$. If in addition $m(\gs) > 0$,
    then the same holds with equality in place of inclusion.
\item
    \label{Item.MContinuity}
    $m(\cdot)$ is uniformly continuous on $\gS_{\chi_0}(m_0)$ for all $m_0 > 0$, $\chi_0 < \infty$.
\item
    \label{Item.MaxContinuity}
	 For all $\gep > 0$, $m_0 > 0$, $\chi_0 < \infty$ there exists $d_0 = d_0 (\chi_0, m_0, \gep) > 0$
    such that if $\gs$, $\wt{\gs} \in \gS_{\chi_0}(m_0)$ satisfy
    $d(\gs, \wt{\gs}) < d_0$ then for all $\gga \in \cL^0$
    there exists $\wt{\gga} \in \cL^0$ such that
    \begin{enumerate}
    \item
        \label{Item.MaxContinuity.1}
        $\normsup{\gga - \wt{\gga}} < \gep$.
    \item
        \label{Item.MaxContinuity.2}
        $\wt{\gs}(\wt{\gga}) > \gs(\gga) - \gep$.
    \item
        \label{Item.MaxContinuity.3}
        $E(\wt{\gga}) < E(\gga) + \gep$.
    \end{enumerate}
\end{enumerate}
\end{prop}
The following are well-known results about large deviation of
simple random walk paths. The emphasis is on the uniformity of the statements. We define $\mu_n \in \cP_n$ as the uniform measure on $\cL^0_n$.
\begin{prop}
\label{Thm.Prop.RWLD}
\forcenewline
\begin{enumerate}
\item
 	\label{Item.UniformPathLDP}
	For any fixed $m_0 > 0$, $\chi_0 < \infty$ as $n \to \infty$ 
	\begin{equation}
 	\label{Eqn.UniformPathLDP}
		-\inv{n} \log \mu_n \lb \gi \subseteq I_{\gs}(s) \rb = E \lb \gGa_{\gs}(\gi) \rb + o(1),
	\end{equation}
	uniformly in all $\gs \in \gS_{\chi_0,n}(m_0)$ and $\gi \in \cI(\gs)$.
\item
 	\label{Item.ConditionalPathLDP}
	For all $\gd > 0$, there exists $\eta = \eta(\gd) > 0$ such that as $n \to \infty$, 
	\begin{equation}
 	\label{Eqn.ConditionalPathLDP}
		-\inv{n} \log \mu_n \lb \normsup{s-\gGa_{\gs}(\gi)} > \gd \rabs \lpr I_{\gs}(s) = \gi \rb
			\geq \eta + o(1),
	\end{equation}
	uniformly in all $\gs \in \gS_{\chi_0,n}(m_0)$ and
	$\gi \in \cI(\gs)$ 
	once $m_0 > 0$ and $\chi_0 < \infty$ are fixed.
\end{enumerate}
\end{prop}

\subsection{Weight-Scaled Environments, Environment Truncation and $\gd$-Optimality}
\label{sub:PreScaled}
For stating the results, it was convenient to work with spatially scaled quantities, such as $\gs_n$ and $\mu_{n, \gb}$. For the proofs, it will turn out useful to define versions of these quantities, which are also weight-scaled. For $\gb \in [0,\infty)$ we set $\ol{\gb} = b_n n^{-1} \gb$. In place of $\gs_n$ and $\mu_{n, \gb}$ we have
\begin{equation}
\label{eqn:mu_def_scaled}
\pi_n = b_n^{-1} \gs_n
\quad , \quad
\ol{\mu}_{n, \gb}(s) =
  \frac{1}{\ol{Q}_{n, \gb}} \, \exp \lb n \gb \pi_n(s) \rb
  \quad \For \ s \in \cL^0_n .
\end{equation}
Clearly $\mu_{n, \gb} = \ol{\mu}_{n, \ol{\gb}}$ and
\eqref{eqn:GbNConditions} reads as 
\begin{equation}
\label{eqn:GBBarAsym}
\ol{\gb}_n \to \ol{\gb}_{\infty} \text{ as } n \to \infty \quad ; \qquad \ol{\gb}_{\infty} \in [0,\infty].
\end{equation}

For $n \leq \infty$, $k < \infty$
we define the {\it truncated} environment
$\pi_n^k$ as the one obtained from $\pi_n$ by removing the masses at indices
$i>k$ (recall that masses are ordered in non-ascending order of their weights). Clearly $\pi_n$, $\pi_n^k$, $\pi_{\infty}$ are random elements of $\gS_{n^2/4-2}$, $\gS_k$, $\gS_{\infty}$ respectively and \eqref{eqn:KMaxsLimit} can be written as
\begin{equation}
\label{eqn:KMaxsLimitAsEnvironments}
\pi_n^k \Ra \pi_{\infty}^k \text{ as } n \to \infty.
\end{equation}
The corresponding polymer measure $\ol{\mu}_{n, \gb}^k$ and the normalization factor
$\ol{Q}_{n, \gb}^k$ are defined as in \eqref{eqn:mu_def_scaled}, but only with
$\pi_n^k$ in place of $\pi_n$.

Next, for $1 \leq k, n \leq \infty$, $\gb \in [0,\infty]$, we define the ``worthiness'' $W_{n, \gb}^k$ of an $\cL^0$-path under environment $\pi_n^k$. If $n = \infty$, we set
\[
    W_{\infty, \gb}^k(\gga) =
        \lbr \begin{array}{ll}
            \gb \pi_{\infty}^k(\gga) - E(\gga)  \quad   & \text{if } \gb < \infty. \\
            \pi_{\infty}^k(\gga)                        & \text{if } \gb = \infty.
        \end{array} \rpr
\]
and if $n < \infty$
\[
    W_{n, \gb}^k(\gga) =
        \lbr \begin{array}{ll}
            \gb \pi_{n}^k(\gga) - E(\gga)
                \quad & \text{if } \ol{\gb}_{\infty} < \infty 
								\text{ or } \gb = 0 .\\
            \pi_{n}^k(\gga) - \inv{\gb} E(\gga)
                        & \text{if } \ol{\gb}_{\infty} = \infty
									\text{ or } \gb = \infty .
        \end{array} \rpr 
\]
We set $\wh{\gga}_{n, \gb}^k$ and $\wh{w}_{n, \gb}^k$ to be the maximizer of $W_{n, \gb}^k$ in $\cL^0$ and its value. It will be shown in Lemma~\ref{lem1} that this is well-defined for all values of $n, \gb, k$, except if $\gb = \infty, n \minwith k<\infty$, in which case a maximizer exist but it is not unique and we apply any a priori deterministic rule for selecting one of the maximizers as $\wh{\gga}_{n, \infty}^k$ (for instance, we may choose the unique one that minimizes $E$). The definition of $\wh{\gga}_{n,\gb}^k$ clearly extends both \eqref{Eqn.FiniteGgaHatDef} and \eqref{Eqn.GgaHatDef} with 
$\wh{\gga}^{\infty}_{n, \ol{\gb}} = \gga_{n, \gb}^*$ and
$\wh{\gga}_{\infty, \gb}^{\infty} = \wh{\gga}_{\infty, \gb}$. From now on we shall omit the superscript $k$ if it is $\infty$. If $\gd > 0$, we also need 
\[
\wh{w}_{n, \gb}^k(\gd) = \maxone{\gga\ \in \cL^0 \,:\; \normsup{\gga - \wh{\gga}_{n, \gb}^k } \geq \gd}
    W_{n, \gb}^k(\gga) .
\]

Finally, we define the {\it remainder} environment as
$\rho_n^k = \pi_n - \pi_n^k$ and its maximal contribution to (minus) the energy of a $\cL^0$ path is $R_n^k = \max_{\gga \in \cL^0} \rho_n^k(\gga)$.

\section{Localization}
\label{sub:localization}
We now prove the main localization results, using 4 lemmas which we state in the beginning of this section. The proofs for these lemmas are deferred to the end of this section, and we first prove Theorems \ref{thm:QuenchedLocalization},\ref{thm:ScalingLimit} and Corollary~\ref{cor:ScalingLimitOfPolymerPath}.

The first lemma establishes the existence and uniqueness of $\wh{\gga}_{n, \gb}^k$ and shows that truncated quantities are good approximations. The first part is due to \cite{Ha07} (Lemma 3.1).
\begin{lem}
\label{lem1}
\forcenewline
\begin{enumerate}
\item
  \label{Item.lem1.0}
	For all $n \leq \infty$ we have 
  $R_n^k \to 0$ as $k \to \infty$ with $\bbP_n$-probability $1$.
\item
  \label{Item.lem1.1}
	For all $1 \leq k, n \leq \infty$, $\gb \in [0,\infty]$ and $\gd > 0$
	quantities  $\wh{\gga}_{n,\gb}^k$, $\wh{w}_{n,\gb}^k$ and 	$\wh{w}_{n,\gb}^k(\gd)$
  are well defined with $\bbP_n$-probability $1$.
  
\item
  \label{Item.lem1.4}
  For all $\gb \in [0,\infty]$, $n \leq \infty$ we have
	$\wh{\gga}_{n,\gb}^k \to \wh{\gga}_{n,\gb} \text{ as } k \to \infty$
  with $\bbP_{n}$-probability $1$.

\item
  \label{Item.lem1.2}
	for all $\gd > 0$, $\gb \in [0,\infty]$ we have
		$\liminf_{k \to \infty} \wh{w}_{\infty,\gb}^k - \wh{w}_{\infty,\gb}^k(\gd) 
		> 0$ with $\bbP_{\infty}$-probability $1$.
\end{enumerate}
\end{lem}

In this lemma we show that truncated quantities of both the finite and limiting system can be coupled such that they are arbitrarily close to each other.
\begin{lem}
\label{Thm.Lem.Coupling}
For all $\gd > 0$, $\gep > 0$, there exist $K$, $(N_k)_{k \geq K}$ such that 
for all $k \geq K$ and $n \geq N_k$ there exists a coupling $\bbP_{n+\infty}$
of $\pi_n^k$ and $\pi_{\infty}^k$ under which with probability at least $1-\gep$:
\begin{enumerate}
\item
    \label{Item.Coupling.1}
    $|\wh{w}_{n,\ol{\gb}_n}^k - \wh{w}_{\infty,\ol{\gb}_{\infty}}^k| < \gep$
\item
    \label{Item.Coupling.2}
    $\normsup{\wh{\gga}_{n,\ol{\gb}_n}^k - \wh{\gga}_{\infty,\ol{\gb}_{\infty}}^k} < \gep$
\item
    \label{Item.Coupling.3}
    $\wh{w}_{n,\ol{\gb}_n}^k(\gd) < \wh{w}_{\infty,\ol{\gb}_{\infty}}^k(\gd/2) + \gep$
\end{enumerate}
\end{lem}

The following improves on the results of Lemma~\ref{lem1} as it shows that approximation by truncated quantities can be made uniform in $n$. The first part is Proposition 3.3 in \cite{Ha07}.
\smallskip
\begin{lem}
\label{lem:UniformFiniteOrderCvgs}
\forcenewline
\begin{enumerate}
\item
\label{item:UnifConvergenceOfRNKInProb}
    For all $\gep, \gd > 0$ there exists $K$ such that $R_n^k < \gd$
    with $\bbP_n$-probability at least $1-\gep$ for all $k > K$ and all $n \leq \infty$.
\item
\label{item:UnifConvergenceOfGGAKInProb}
    For all $\gep,\gd > 0$ there exists $K$ such that $\normsup{\wh{\gga}_{n,\ol{\gb}_n}^k - \wh{\gga}_{n,\ol{\gb}_n}} < \gd$ with $\bbP_n$-probability at least $1-\gep$ for all $k > K$ and all $n \leq \infty$.
\item
\label{item:UnifGDComparison}
    For all $\gep, \gd > 0$, there exists $K, \eta > 0$ such that
    $\wh{w}_{n, \ol{\gb}_n}^k(\gd) < \wh{w}_{n, \ol{\gb}_n}^k - \eta$
    with $\bbP_n$-probability at least $1-\gep$
    for all $k > K$ and $k \leq n \leq \infty$.
\end{enumerate}
\end{lem}

In this lemma we show concentration with truncated quantities (unless  $\ol{\gb}_{\infty} = \infty$, in which case this is essentially
Theorem~\ref{thm:QuenchedLocalization}).
\begin{lem}
\label{lem:ConcentrationOnTruncated}
For all $\gd > 0$, $\gep >0$, there exists $\nu > 0$ and $K, (N_k)_{k \geq K}$ (case $\ol{\gb}_{\infty} < \infty$) or $N$ (case $\ol{\gb}_{\infty} = \infty$)
such that with $\bbP_n$-probability at least $1-\gep$
\begin{equation}
\label{Eqn.Lemma2Proof.6}
\begin{array}{ll}
-\inv{n} \log \ol{\mu}_{n, \ol{\gb}_n}^k \lb s\,:\; \normsup{ s - \wh{\gga}_{n, \ol{\gb}_n}^k } > \gd \rb \geq \nu
    \qquad \qquad & \text{(case } \ol{\gb}_{\infty} < \infty \text{)} \\
-\inv{n} \log \ol{\mu}_{n, \ol{\gb}_n} \lb s\,:\; \normsup{ s - \wh{\gga}_{n, \ol{\gb}_n} } > \gd \rb \geq \nu \ol{\gb}_n
    & \text{(case } \ol{\gb}_{\infty} = \infty \text{)},
\end{array}
\end{equation}
for all $k > K$, $n > N_k$ (case $\ol{\gb}_{\infty} < \infty$) or $n > N$ (case $\ol{\gb}_{\infty} = \infty$).
\end{lem}

\medskip
\begin{proof}[Proof of Theorem ~\ref{thm:QuenchedLocalization}]
If $\ol{\gb}_{\infty} = \infty$ we can just quote Lemma~\ref{lem:ConcentrationOnTruncated}. Otherwise, fix $\gep>0$, $\gd>0$ and write
\begin{eqnarray}
\nonumber
\ol{\mu}_{n, \ol{\gb}_n} \lb s\,:\; \normsup{ s - \wh{\gga}_{n, \ol{\gb}_n}^k } > \gd \rb
	& \leq  &
		\ol{\mu}_{n, \ol{\gb}_n}^k \lb s\,:\; \normsup{ s - \wh{\gga}_{n, \ol{\gb}_n}^k } > \gd \rb \\
\label{eqn:Lem3Derivative}
	&       & \quad
		\times \ \sup \lbr \frac{d \ol{\mu}_{n, \ol{\gb}_n}}{d \ol{\mu}_{n, \ol{\gb}_n}^k}(s) \;:\;\;
							\normsup{ s - \wh{\gga}_{n, \ol{\gb}_n}^k } > \gd \rbr
\end{eqnarray}
By Lemma~\ref{lem:ConcentrationOnTruncated} the first factor on the r.h.s. is exponentially decaying in $n$ with some rate $\nu > 0$ with probability at least $1-\gep$ for all properly large $k$, $n$. On the other hand
\begin{equation*}
\frac{d \ol{\mu}_{n, \ol{\gb}_n}}{d \ol{\mu}_{n, \ol{\gb}_n}^k}(s)
	\leq  \frac{\ol{Q}_{n, \ol{\gb}_n}^k}{\ol{Q}_{n, \ol{\gb}_n}} \mbox{exp} \lb \ol{\gb}_n n R_n^k \rb \\
	\leq  \mbox{exp} \lb \ol{\gb}_n n R_n^k \rb,
\end{equation*}
and therefore using Lemma~\ref{lem:UniformFiniteOrderCvgs} part \eqref{item:UnifConvergenceOfRNKInProb} and choosing $k$ large enough, we can have the second factor in \eqref{eqn:Lem3Derivative} grow exponentially in $n$ with rate at most $\nu/2$ also with probability at least $1-\gep$.
Finally, from part \eqref{item:UnifConvergenceOfGGAKInProb} of Lemma~\ref{lem:UniformFiniteOrderCvgs} for possibly larger $k$, we can have also $\normsup{\wh{\gga}_{n, \ol{\gb}_n} - \wh{\gga}_{n, \ol{\gb}_n}^k} < \gd/2$ with the same probability. Combining the above, we complete the proof.
\end{proof}

\smallskip
\begin{proof}[Proof of Theorem~\ref{thm:ScalingLimit}]
Fix $\gep, \gd > 0$. From Lemma \ref{lem1}, ~\ref{Thm.Lem.Coupling} and \ref{lem:UniformFiniteOrderCvgs} it follows that we can find $k$ large enough and then $n$ large enough such that with $\bbP_{n+\infty}$ probability at least $1-\gep$
\[
\normsup{\wh{\gga}_{n, \ol{\gb}_n} - \wh{\gga}^k_{n, \ol{\gb}_n}} < \gd
\,, \quad
\normsup{\wh{\gga}^k_{n, \ol{\gb}_n} - \wh{\gga}^k_{\infty, \ol{\gb}_{\infty}}} < \gd
\,, \quad
\normsup{\wh{\gga}^k_{\infty, \ol{\gb}_{\infty}} - \wh{\gga}_{\infty, \ol{\gb}_{\infty}}} < \gd
\,.
\]
This gives $\normsup{\wh{\gga}_{n, \ol{\gb}_n} - \wh{\gga}_{\infty, \ol{\gb}_{\infty}}} < 3\gd$ with $\bbP_{n+\infty}$ probability $1-\gep$
for all $n$ sufficiently large. Since $\gep$, $\gd$ are arbitrary, the result follows.
\end{proof}

\smallskip
\begin{proof}[Proof of Corollary~\ref{cor:ScalingLimitOfPolymerPath}]
We can use Skorohod Representation Theorem or the proof of Theorem~\ref{thm:ScalingLimit}, to conclude that for any $\gep, \gd > 0$ if $n$ is large enough
$\bbP_{n+\infty} (\normsup{\wh{\gga}_{n, \ol{\gb}_n} - \wh{\gga}_{\infty, \ol{\gb}_{\infty}}} < \gd) > 1-\gep$. Then, possibly for larger $n$, from Theorem~\ref{thm:QuenchedLocalization} we have
$\ol{\mu}_{n, \ol{\gb}_n} \lb \normsup{s-\wh{\gga}_{n, \ol{\gb}_n}} > \gd \rb \leq \gep$ with
$\bbP_n$-probability at least $1-\gep$. By the total probability formula this implies $\bbP_n(\normsup{S_{n, \gb_n} - \wh{\gga}_{n, \ol{\gb}_n}} > \gd) \leq 2\gep$. All together we have
$\bbP_{n+\infty} \lb \normsup{ S_{n, \gb_n} - \wh{\gga}_{\infty, \ol{\gb}_{\infty}} } > 2\gd \rb \leq 3\gep$ 
and since $\gd, \gep$ are arbitrary, the result follows.
\end{proof}

\smallskip
\begin{proof}[Proof of Lemma ~\ref{lem1}]
Part \eqref{Item.lem1.0} is Lemma 3.1 in \cite{Ha07}.
For part \eqref{Item.lem1.1}, existence of a maximizer in \eqref{Eqn.GgaHatDef} will follow if we show that
$W_{n,{\gb}}^k$ is upper semi-continuous, since $\cL^0$ is compact in the $\normsup{\cdot}$ topology. Indeed, $E(\cdot)$ is lower semi-continuous (Proposition~\ref{prop:EProps}). As for $\pi_n^k(\cdot)$,
from part \eqref{Item.lem1.0}, given $\gep > 0$ and $\gga_0 \in \cL^0$ we may find $k_1$ such that
$R_n^{k_1} < \gep$, and set $\gd = \min \lbr |y^l_{\gs} - \gga_0(x^l_{\gs})| : 1 \leq l < k_1 ,\, y^l_{\gs} \neq \gga_0(x^l_{\gs}) \rbr$.
Then, for $\normsup{\gga - \gga_0} < \gd$
\begin{eqnarray*}
\pi_n^k(\gga) < \pi_n^k(\gga_0) + R_n^{k_1} < \pi_n^k(\gga_0) + \gep
\end{eqnarray*}
which implies upper semi-continuity.

It remains to show uniqueness.
If $\gb = 0$ we have $\wh{\gga}_{n, 0}^k \equiv 0$, which is the minimum of $E$ (Proposition~\ref{prop:EProps}). If $\gb = n = k = \infty$, uniqueness was proved in \cite{Ha07} (see Proposition 4.1 and the unique way to extend $U^*$ to a continuous increasing path $P^*$). For $\gb=\infty$, $n \minwith k < \infty$ uniqueness holds by definition. In the remaining cases, assume the contrary and let $\wh{\gga}_0$, $\wc{\gga}_0$
be two different maximizers of $W_{n, \gb}^k$. Then $\wh{\gi}_0 = I_{\pi_n^k}(\wh{\gga_0})$
must be different from $\wc{\gi}_0 = I_{\pi_n^k}(\wc{\gga_0})$, because the minimizer of $E(\gga)$ in $\{\gga \in \cL^0 : I_{\pi_n^k}(\gga) = \gi\}$ is unique, as it follows
from the strict convexity of $E(\cdot)$.
We proceed as in Proposition 4.1 in \cite{Ha07}.
Without loss of generality there must exists $1 \leq j \leq n$ such that with positive $\bbP_n$-probability
\[
\max_{\gga\, :\; j \in I_{\pi_n^k}(\gga)} W_{n, \gb}^k (\gga) =
\max_{\gga\, :\; j \notin I_{\pi_n^k}(\gga)} W_{n, \gb}^k (\gga) .
\]
But conditioned on $(Z_n^i)_{i \geq 1}$, $(U_n^i)_{i \neq j}$
the r.h.s of the above event is a constant while the l.h.s is an absolutely continuous
(w.r.t. Lebesgue measure) random variable. It follows that this probability is zero,
which is a contradiction. This proves part \eqref{Item.lem1.1}.

Part \eqref{Item.lem1.4} is is trivial if $n < \infty$. If $n=\infty$, by compactness
$\exists (k_l)_{l \geq 1}$, $\wt{\gga}_0 \in \cL^0$ s.t. $\wh{\gga}_{\infty, \gb}^{k_l} \to \wt{\gga}_0$
as $l \to \infty$.
Then
\begin{eqnarray*}
  W_{\infty, \gb} \lb \wt{\gga}_0 \rb
	& \geq	&	\limsup_{l \to \infty} W_{\infty, \gb} \lb \wh{\gga}_{\infty, \gb}^{k_l} \rb
					\geq \limsup_{l \to \infty} W_{\infty, \gb}^{k_l} \lb \wh{\gga}_{\infty, \gb}^{k_l} \rb \\
  	& \geq  &	\limsup_{l \to \infty} W_{\infty, \gb}^{k_l} \lb \wh{\gga}_{\infty, \gb} \rb
					= W_{\infty, \gb} \lb \wh{\gga}_{\infty, \gb} \rb
					= \wh{w}_{\infty, \gb},
   \end{eqnarray*}
where the first inequality follows from upper semi-continuity.
By uniqueness it must be that $\wt{\gga}_0 = \wh{\gga}_{\infty, \gb}$ and since this is true
for any subsequence of $\wh{\gga}_{\infty, \gb}^k$, the result follows.

As for part \eqref{Item.lem1.2}, if the statement had been false, then there would have been sequences $(k_l)_{l \geq 1}$ and $(\wc{\gga}^{k_l}_{\infty, \gb})_{l \geq 1}$ such that
$\limsup_{l \to \infty} W_{\infty, \gb}^{k_l} (\wc{\gga}^{k_l}_{\infty, \gb}) \geq \liminf_{l \to \infty} W_{\infty, \gb}^{k_l}(\wh{\gga}^{k_l}_{\infty, \gb})$ and $\normsup{\wc{\gga}^{k_l}_{\infty, \gb} - \wh{\gga}_{\infty, \gb}^{k_l}} \geq \gd$.
By compactness we could further suppose that $\wc{\gga}^{k_l}_{\infty, \gb} \to \wc{\gga}_{\infty, \gb}$ as $l \to \infty$. Then upper semi-continuity would have implied
\begin{eqnarray*}
	W_{\infty, \gb}(\wc{\gga}_{\infty, \gb}) 
		& \geq & \limsup_{l \to \infty} W_{\infty, \gb}(\wc{\gga}^{k_l}_{\infty, \gb})
		\geq \limsup_{l \to \infty} W_{\infty, \gb}^{k_l}(\wc{\gga}^{k_l}_{\infty, \gb}) \\
	& \geq & \liminf_{l \to \infty} W_{\infty, \gb}^{k_l}(\wh{\gga}^{k_l}_{\infty, \gb}) 
	\geq \liminf_{l \to \infty} W_{\infty, \gb}^{k_l}(\wh{\gga}_{\infty, \gb}) 
	= W_{\infty, \gb}(\wh{\gga}_{\infty, \gb}) ,
\end{eqnarray*}
and part \eqref{Item.lem1.4} would have given  $\normsup{\wc{\gga}_{\infty, \gb} - \wh{\gga}_{\infty, \gb}} \geq \gd > 0$. This would have violated the uniqueness of the global maximizer.
\end{proof}

\smallskip
\begin{proof}[Proof of Lemma~\ref{Thm.Lem.Coupling}]
As in Proposition 3.2 of \cite{Ha07}, it follows from \eqref{eqn:KMaxsLimitAsEnvironments} and Skorohod Representation Theorem that for any $k, d_0 > 0$, we can couple together $\pi_n^k$ and $\pi_{\infty}^k$ such that $d(\pi_n^k, \pi_{\infty}^k) < d_0$ with arbitrarily high probability as long as $n$ is large enough. Call this coupling measure $\bbP_{n+\infty}$ and observe that
the absolute continuity of $Z^j_{\infty}, j=1, \dots k$ implies that by choosing $m_0$ sufficiently small, we can make $m(\pi_{\infty}^k) > m_0$ occur with $\bbP_{n+\infty}$ probability arbitrarily close to $1$. Using Proposition~\ref{Thm.Prop.Continuity} part \eqref{Item.MContinuity}, if $d_0$ is chosen small enough, this implies
$\pi_n^k, \pi_{\infty}^k \in \gS_k(m_0/2)$.
Now given $\gep>0$, Proposition~\ref{Thm.Prop.Continuity} part \eqref{Item.MaxContinuity} applied to $\pi_n^k$ and $\pi_{\infty}^k$ with  $\wh{\gga}_{n, \ol{\gb}_n}^k$ and $\wh{\gga}_{\infty, \ol{\gb}_{\infty}}^k$  together with the assumption on $\ol{\gb}_n$ guarantee that by further restricting $d_0$ we have $|\wh{w}_{n, \ol{\gb}_n}^k - \wh{w}_{\infty, \ol{\gb}_{\infty}}^k | < \gep$. This shows \eqref{Item.Coupling.1}.

As for \eqref{Item.Coupling.2}, for any $\gep > 0$ it follows from Lemma~\ref{lem1} part \eqref{Item.lem1.2}, that there exists $\eta > 0$, such that $\wh{w}_{\infty, \ol{\gb}_{\infty}}^k(\gep) < \wh{w}_{\infty, \ol{\gb}_{\infty}}^k - \eta$, if $k$ is large enough with arbitrarily high $\bbP_{n+\infty}$ probability.
Then, using Proposition ~\ref{Thm.Prop.Continuity} part \eqref{Item.MaxContinuity} again, for
$d_0$ and $n$ large it must be that $W_{n, \ol{\gb}_n} ^k(\gga) < \wh{w}_{\infty, \ol{\gb}_{\infty}}^k - \eta/2$  for all $\gga \in \cL^0$ such that $\normsup{\gga - \wh{\gga}_{\infty, \ol{\gb}_{\infty}}^k} > 2 \gep $. Then from part \eqref{Item.Coupling.1}, for possibly smaller $d_0$, all such $\gga$ satisfy $W_{n, \ol{\gb}_n}^k(\gga) < \wh{w}_{n, \ol{\gb}_n}^k - \eta/3$. This shows $\normsup{\wh{\gga}_{n, \ol{\gb}_n}^k -\wh{\gga}_{\infty, \ol{\gb}_{\infty}}^k} \leq 2 \gep$ as required. 

Finally, given $\gd > 0$ we may find $\gga \in \cL^0$ satisfying
$\normsup{\gga - \wh{\gga}_{n, \ol{\gb}_n}^k} \geq \gd$ and $W_{n, \ol{\gb}_n}^k(\gga) =
\wh{w}_{n, \ol{\gb}_n}^k(\gd)$. Then using part \eqref{Item.Coupling.2} and Proposition~\ref{Thm.Prop.Continuity} part \eqref{Item.MaxContinuity}, provided $d_0$ is sufficiently small and $k,n$ are sufficiently large, we may find $\wt{\gga} \in \cL^0$ satisfying
$\normsup{\wt{\gga} - \wh{\gga}_{\infty, \ol{\gb}_{\infty}}^k} \geq \gd/2$ and $W_{\infty, \ol{\gb}_{\infty}}^k(\wt{\gga}) >
\wh{w}_{n, \ol{\gb}_n}^k(\gd) - \gep$ where $\gep > 0$ is given. This shows \eqref{Item.Coupling.3}.
\end{proof}

\begin{proof}[Proof of Lemma~\ref{lem:UniformFiniteOrderCvgs}]
Part \eqref{item:UnifConvergenceOfRNKInProb} is Proposition 3.3 in \cite{Ha07}.
For part \eqref{item:UnifConvergenceOfGGAKInProb}, fix $\gd, \gep > 0$ and
use Lemma~\ref{lem1} part \eqref{Item.lem1.2}
to find $\eta > 0$ such that
$\wh{w}_{\infty, \ol{\gb}_{\infty}}^k(\gd) < \wh{w}_{\infty, \ol{\gb}_{\infty}}^k - \eta$ with $\bbP_{\infty}$-probability at least  $1-\gep$, as long as $k$ is large enough. Then, use Lemma \ref{Thm.Lem.Coupling}, to find $n$ large enough such that
\begin{equation}
\label{eqn:A1}
\wh{w}_{n, \ol{\gb}_n}^k(2\gd) < \wh{w}_{n, \ol{\gb}_n}^k - \eta/2 \,
\end{equation}
with $\bbP_{n+\infty}$-probability at least $1 - 2\gep$.
Now, from part \eqref{item:UnifConvergenceOfRNKInProb},
with probability at least $1-3\gep$
we can also have $R_n^k < \eta/4$ for all $n$, possibly by further restricting $k$. In this case, it must be that
\begin{equation}
\label{eqn:A2}
\normsup{\wh{\gga}_{n, \ol{\gb}_n}^l - \wh{\gga}_{n, \ol{\gb}_n}^k} < 2\gd \,,
\end{equation}
for all $l \geq k$. Then
$\normsup{\wh{\gga}_{n, \ol{\gb}_n}^l - \wh{\gga}_{n, \ol{\gb}_n}} < 4\gd$ for all $l \geq k$. This shows part \eqref{item:UnifConvergenceOfGGAKInProb}. But then from \eqref{eqn:A1}, \eqref{eqn:A2} and the restriction on $R_n^k$, we have $\wh{w}_{n, \ol{\gb}_n}^l(4\gd) < \wh{w}_{n, \ol{\gb}_n}^l - \eta/4$ for all $l \geq k$, for possibly larger $k$. This shows part \eqref{item:UnifGDComparison}. The restriction on $n$ can be enforced by a restriction on $k$. 
\end{proof}

\begin{proof}[Proof of Lemma~\ref{lem:ConcentrationOnTruncated}]
Let us treat the $\ol{\gb}_{\infty} < \infty$ case first. Write
\begin{eqnarray}
\nonumber
\lefteqn{
\ol{\mu}_{n,\ol{\gb}_n}^k ( s\,:\; \normsup{ s - \wh{\gga}_{n, \ol{\gb}_n}^k } > \gd )} \\
\nonumber
	& \leq  &
		\ol{\mu}_{n,\ol{\gb}_n}^k ( s\,:\; \normsup{ \gL_{\pi_n^k}(s) - \wh{\gga}_{n, \ol{\gb}_n}^k } \geq \gd/2 ) \\
\label{eqn:Lem3Splitting}
	&       &
		\; + \; \ol{\mu}_{n,\ol{\gb}_n}^k ( s\,:\; \normsup{ s - \gL_{\pi_n^k}(s) } > \gd/2\;,
			\;\; \normsup{ \gL_{\pi_n^k}(s) - \wh{\gga}_{n, \ol{\gb}_n}^k } < \gd/2 ).
\end{eqnarray}
Now,
\begin{equation}
\label{Eqn.Proposition2Proof.3}
\ol{\mu}_{n,\ol{\gb}_n}^k \lb s\,:\; \normsup{ \gL_{\pi_n^k}(s) - \wh{\gga}_{n, \ol{\gb}_n}^k } \geq \gd/2 \rb
	\leq \sum_{\gi}
	\ol{\mu}_{n,\ol{\gb}_n}^k \lb I_{\pi_n^k}(s) = \gi \rb	
\end{equation}
where the sum is over all $\gi \in \cI(\pi_n^k)$ such that $\normsup{\gGa_{\pi_n^k}(\gi) - \wh{\gga}_{n, \ol{\gb}_n}^k} \geq \gd/2$. Each
such term satisfies
\begin{eqnarray}
\nonumber
\ol{\mu}_{n,\ol{\gb}_n}^k \lb I_{\pi_n^k}(s) = \gi \rb
	& = 	& \frac{\mu_n \lsp \exp \lb n \ol{\gb}_n \pi_n^k(s) \rb ;\;I_{\pi_n^k}(s) = \gi \rsp}
					   {\mu_n \lsp \exp \lb n \ol{\gb}_n \pi_n^k(s) \rb \rsp}		 \\
\label{eqn:OnePathContrib}
	& \leq	& \frac	{\exp \lb n \ol{\gb}_n \pi_n^k \lb \gGa_{\pi_n^k} (\gi) \rb \rb
						\mu_n \lb \gi \subseteq I_{\pi_n^k}(s) \rb}
					{\exp \lb n \ol{\gb}_n \pi_n^k \lb \wh{\gga}_{n, \ol{\gb}_n}^k \rb \rb
						\mu_n \lb I_{\pi_n^k} \lb \wh{\gga}_{n, \ol{\gb}_n}^k \rb \subseteq I_{\pi_n^k}(s) \rb}	.
\end{eqnarray}
As in the previous proof, for any $k>0$ we may find $m_0 > 0$, such that with arbitrarily high probability
$\pi_n^k \in \gS_k(m_0)$ as long as $n$ is large enough.
Then, by Proposition~\ref{prop:EProps} part \eqref{item:EProps2} and
Proposition \ref{Thm.Prop.RWLD} part \eqref{Item.UniformPathLDP} 
\begin{eqnarray*}
-\inv{n} \log \lsp \ol{\mu}_{n,\ol{\gb}_n}^k \lb I_{\pi_n^k}(s) = \gi \rb \rsp 	
	& \geq 	& W_{n, \ol{\gb}_n}^k \lb \gL_{\pi_n^k} \lb \wh{\gga}_{n, \ol{\gb}_n}^k \rb
				\rb - W_{n, \ol{\gb}_n}^k \lb \gGa_{\pi_n^k} (\gi) \rb + o(1)	 \\
	& \geq 	& \wh{w}_{n, \ol{\gb}_n}^k - W_{n, \ol{\gb}_n}^k \lb \gGa_{\pi_n^k} (\gi) \rb + o(1).
\end{eqnarray*}
Plugging this into \eqref{Eqn.Proposition2Proof.3}, noting that
there are at most $2^k$ terms in the sum there and using Lemma~\ref{lem:UniformFiniteOrderCvgs} part \eqref{item:UnifGDComparison},
we infer that there exists $\nu_1 > 0$, such that
\begin{equation}
\label{Eqn.Proposition2Proof.5}
\bbP_n \lb -\inv{n} \log \ol{\mu}_{n,\ol{\gb}_n}^k \lb s\,:\; \normsup{ \gL_{\pi_n^k}(s) - \wh{\gga}_{n, \ol{\gb}_n}^k } \geq \gd/2 \rb
	\geq \nu_1 \rb \geq 1-\gep
\end{equation}
for large enough $k$, $n$.
At the same time, the second term in (\ref{eqn:Lem3Splitting}) is clearly bounded by
\[
	\sum_{\gi}
			\ol{\mu}_{n,\ol{\gb}_n}^k \lb \lpr \normsup{ s - \gGa_{\pi_n^k}(\gi) } > \gd/2 \rabs
				I_{\pi_n^k}(s) = \gi \rb 
=\sum_{\gi}
			\mu_n \lb \lpr \normsup{ s - \gGa_{\pi_n^k}(\gi) } > \gd/2 \rabs
				I_{\pi_n^k}(s) = \gi \rb ,
\] 
where the sum is over all $\gi \in \cI(\pi_n^k)$ such that $\normsup{\gGa_{\pi_n^k}(\gi) - \wh{\gga}_{n, \ol{\gb}_n}^k} < \gd/2$. 
We may now use Proposition~\ref{Thm.Prop.RWLD} part \eqref{Item.ConditionalPathLDP} and the bound on the number of terms, to conclude that there exists $\nu_2 > 0$ such that
\begin{equation}
\label{Eqn.Proposition2Proof.6}
\bbP_n \lb	
-\inv{n} \log \lsp \ol{\mu}_{n,\ol{\gb}_n}^k \lb s\,:\; \normsup{ s - \gL_{\pi_n^k}(s) } > \gd/2\;,
		\;\; \normsup{ \gL_{\pi_n^k}(s) - \wh{\gga}_{n, \ol{\gb}_n}^k } < \gd/2 \rb \rsp 	
		\geq \nu_2 \rb \geq 1-\gep
\end{equation}
This holds uniformly in $k$, but $n$ needs to be large enough. Combining \eqref{eqn:Lem3Splitting},
\eqref{Eqn.Proposition2Proof.5} and \eqref{Eqn.Proposition2Proof.6}, we complete the $\ol{\gb}_{\infty} < \infty$ case.

If $\ol{\gb}_{\infty} = \infty$ we take $k=|\cD_n^0| = n^2/4-2$. Then from \eqref{eqn:OnePathContrib} we get
for all $\gga \in \cL^0_n$
\[
-\inv{n} \log \ol{\mu}_{n,\ol{\gb}_n} (\gga)
    \geq \ol{\gb}_n (\wh{w}_{n, \ol{\gb}_n} - W_{n, \ol{\gb}_n} (\gga)) + O(1) \,.
\]
Then, since there are at most $2^n$ paths and using Lemma~\ref{lem:UniformFiniteOrderCvgs} part \eqref{item:UnifGDComparison}, we obtain for some $\nu > 0$
\[
-\inv{n} \log \ol{\mu}_{n,\ol{\gb}_n} \lb s\,:\; \normsup{ s - \wh{\gga}_{n, \ol{\gb}_n} } > \gd \rb
    \geq \nu \ol{\gb}_n - O(1) \geq \tfrac{1}{2} \nu \ol{\gb}_n \,,
\]
with $\bbP_n$-probability at least $1-\gep$ as long as $n$ is large enough.
\end{proof}

\section{The Limit distribution}
\label{Section:PropLimitDist}
\begin{proof}[Proof of Proposition~\ref{prop:MIsDistingushiable}]
For $\gb \in (0,\infty)$, set $Y_\gb(\gga) = \gb^{-1} W_{\infty, \gb}(\gga) = \pi_{\infty}(\gga) - \gb^{-1} E(\gga)$. If $0 <\gb_2 < \gb_1 < \infty$ then
\[
E(\wh{\gga}_{\infty, \gb_2}) =
    \frac{Y_{\gb_1}(\wh{\gga}_{\infty, \gb_2}) - Y_{\gb_2}(\wh{\gga}_{\infty, \gb_2})}
        {\gb_2^{-1} - \gb_1^{-1}} \leq
    \frac{Y_{\gb_1}(\wh{\gga}_{\infty, \gb_1}) - Y_{\gb_2}(\wh{\gga}_{\infty, \gb_1})}
        {\gb_2^{-1} - \gb_1^{-1}} = E(\wh{\gga}_{\infty, \gb_1})\,.
\]
This shows $E(\wh{\gga}_{\infty, \gb_2}) \le_s E(\wh{\gga}_{\infty, \gb_1})$. It remains to show that $E(\wh{\gga}_{\infty, \gb_2}) < E(\wh{\gga}_{\infty, \gb_1})$ with positive probability.
Recall the definitions of $Z^i$, $V^i$ and $T_i$ in \eqref{eqn:ExtremesLimit}.
For any $\gd > 0$, from Lemma~\ref{lem1} part \eqref{Item.lem1.0} we may find $k_0$ such that $R_{\infty}^{k_0} < \gd$ with $\bbP_{\infty}$-probability at least $\inv{2}$. Also, there exists $t_0$ large enough such that $T_{k_0} < t_0$ with $\bbP_{\infty}$-probability at least $\tfrac{3}{4}$. Then
\[
\inv{2} \leq \bbP_{\infty} (R_{\infty}^{k_0} < \gd) =
    \bbE_{\infty} \lb \bbP_{\infty} (R_{\infty}^{k_0} < \gd | T_{k_0}) \rb
        \leq \bbP_{\infty} (R_{\infty}^{k_0} < \gd | T_{k_0} = t_0) + \inv{4}
\]
where for the last inequality we use the obvious monotonicity of $t \mapsto \bbP_{\infty} (R_{\infty}^{k_0} < \gd | T_{k_0} = t)$. 
Since the last term is independent of $k_0$ and using the monotonicity again we in fact have
\begin{equation}
\label{eqn:SmallRWithPosProb}
\bbP_{\infty} (R_{\infty}^{k} < \gd | T_{k} \geq t_0) \geq \inv{4}
\end{equation}
for all $k$. Next, it is not difficult to verify that we can find $\gd, \eta > 0$ small enough as well as $0 < t_1 < t_2 < t_3$ and $z_1, z_2 \in \cD^0$, such that on the event
\[
\cA = \lbr |T_i - t_i| < \eta, \norm{Z^i - z_i} < \eta\text{ for } i=1,2
    \ ,\, T_3 \geq t_3 \text{ and } R^3_{\infty} < \gd \rbr
\]
the following holds:
\begin{enumerate}
\item
    $|a_{\pi_{\infty}}(1, 2)| > 1$. In other words, no $\cL^0$-curve can take both $Z^1$ and $Z^2$.
\item
    $E(\gga_1) > E(\gga_2)$ and $W_{\infty, \gb_i}(\gga_j) > 0$ for $i,j = 1,2$,
	where $\gga_i = \gGa_{\pi_{\infty}}(\{0, i, \infty\})$.
\item
	$\exists \gep > 0$ such that
    $W_{\infty, \gb_1}(\gga_1) - W_{\infty, \gb_1}(\gga_2) > \gep$ and
    $W_{\infty, \gb_2}(\gga_2) - W_{\infty, \gb_2}(\gga_1) > \gep$.
\item
    $\gb_1(V^3 + R^3_{\infty}) < \gep$.
\end{enumerate}
In light of \eqref{eqn:ExtremesLimit} and \eqref{eqn:SmallRWithPosProb}, event $\cA$ has positive probability under $\bbP_{\infty}$ for any choice of parameters. At the same time, the above conditions guarantee that $\wh{\gga}_{\infty, \gb_1} \neq \wh{\gga}_{\infty, \gb_2}$. Now,
if $E(\wh{\gga}_{\infty, \gb_1}) = E(\wh{\gga}_{\infty, \gb_2})$ it must be that $\pi_{\infty}(\wh{\gga}_{\infty, \gb_1}) = \pi_{\infty}(\wh{\gga}_{\infty, \gb_2})$, for otherwise one of $\wh{\gga}_{\infty, \gb_i}$ cannot be a maximizer. But then it has to be the case that $\wh{\gga}_{\infty, \gb_1}, \wh{\gga}_{\infty, \gb_2}$ are both maximizers of $W_{\infty, \gb_1}$ and $W_{\infty, \gb_2}$ which contradicts uniqueness. Therefore
$E(\wh{\gga}_{\infty, \gb_1}) > E(\wh{\gga}_{\infty, \gb_2})$ as desired. The cases $\gb_2=0$ and/or $\gb_1=\infty$ are proved in a similar way.
\end{proof}

\begin{proof}[Proof of Proposition~\ref{prop:PhaseTransitionAtZero}]
Although $\wh{\gga}_{\infty,\gb}$ might not be defined for all $\gb$ for a given environment, it is the case for $\wh{w}_{\infty,\gb}$ with $\bbP_{\infty}$-probability $1$, as the proof of Lemma~\ref{lem1} shows. This makes $\gb_c$ well defined and shows part 
\eqref{item:PhaseTransitionAtZero0}.
Part \eqref{item:PhaseTransitionAtZero1} holds since if for some $\gb_0$ and $\gga_0 \in \cL^0 \setminus \{0\}$ we have $W_{\infty, \gb_0}(\gga_0) \geq 0$ then $\wh{w}_{\infty, \gb} > 0$ for all $\gb > \gb_0$. For parts \eqref{item:PhaseTransitionAtZero2} and
\eqref{item:PhaseTransitionAtZero3}, define $\gga_{z} \in \cL^0$ as the curve $\gga_z = \linear \{(0,0), z, (1,0)\}$, where $z=(x,y) \in \cD$. We claim that
\[
C_1 \lb \tfrac{y^2}{x} + \tfrac{y^2}{1-x} \rb \leq E(\gga_z) \leq C_2 \lb \tfrac{y^2}{x} + \tfrac{y^2}{1-x} \rb
\]
for some positive $C_1, C_2$. This can be verified by a simple calculation.
This in turn implies that the set
$\{z \in \cD :\:
E(\gga_z) \leq \gd \}$ has Lebesgue measure $\gT(\sqrt{\gd})$ as $\gd \to 0$.

Next, from LLN there exists a.s. $k_0$ such that $T_k \leq 2k$ for all $k > k_0$. Then, conditioning on $(T_k)_{k \geq 1}$, for all such $k \geq k_0$ and any $\gb >0$ we have
\[
\bbP_{\infty} \lb \lpr W_{\infty, \gb}(\gga_{Z^k}) > 0 \rabs (T_k)_{k\geq 1} \rb \geq
\bbP_{\infty} \lb E (\gga_{Z^k}) < \gb (2k)^{-\inv{\ga}} \rb \geq
C \gb^{1/2} k^{-\inv{2\ga}} 
\]
and events $\{ W_{\infty, \gb}(\gga_{Z^k}) > 0 \}_{k \geq k_0}$ are (conditionally) independent. Now, for $\ga > \inv{2}$ and any $\gb > 0$, the sum
of the probabilities above diverges, whence we may conclude via Borel-Cantelli Lemma that with $\bbP_{\infty}$-probability $1$ there will be $k_1 > k_0$ for which $0 < W_{\infty, \gb}(\gga_{Z^{k_1}}) \leq \wh{w}_{\infty, \gb}$. Since $\gb$ is arbitrary the proof for part \eqref{item:PhaseTransitionAtZero2} is complete.

For part \eqref{item:PhaseTransitionAtZero3}, we need to show $\bbP_{\infty} (\wh{w}_{\infty, \gb} = 0) \to 1$ as $\gb \to 0$. Let $\gep > 0$ be arbitrarily small. It is not difficult to see that there exists $\eta > 0$ such that $T_k \geq \eta k$ for all $k \geq 1$ with $\bbP_{\infty}$-probability at least $1-\gep$. On this event and if $\ga < 1$
\[
Q^k \bydef \sum_{m \geq k} V^m \leq C k^{-\lb\inv{\ga} - 1\rb}.
\]
Now, for $\gga \in \cL^0$ set $l(\gga) = \min \lb I_{\pi_{\infty}}(\gga) \setminus \{0\} \rb$ - the smallest index of a mass reached by $\gga$. Then Proposition~\ref{prop:EProps}-\eqref{item:EProps2} implies that $W_{\infty, \gb} (\gga) \leq \gb Q^{l(\gga)} - E(\gga_{Z^{l(\gga)}})$ 
for all $\gga \in \cL^0$. Therefore, if $\ga$ is further restricted $\ga < \inv{3}$ and $\gb$ is small enough,
\begin{eqnarray*}
\bbP_{\infty} (\wh{w}_{\infty, \gb} = 0) +\gep
    & \geq  & \bbP_{\infty} (\gb Q^k - E(\gga_{Z^k}) \leq 0 \text{ for all } k
		\ |\  T_k \geq \eta k \ ; \; k = 1,\dots )  \\
    & \geq  & \prod_{k=1}^{\infty} \lb 1 - C \gb^{1/2}
        k^{-\lb\inv{2\ga} - \inv{2}\rb} \rb \\
    & \geq  & \exp \lbr -C \gb^{1/2} \sum_{k=1}^{\infty} k^{-\lb\inv{2\ga} - \inv{2}\rb} \rbr
\end{eqnarray*}
and the last term goes to $1$ as $\gb \to 0$. Since $\gep$ is arbitrarily, this concludes the last part of the proposition. 
\end{proof}

\section{Proofs for Subsection~\ref{Outline.Sub.EnvironmentsAndLipCurves}}
\label{sub:DeferredProofs}
\begin{proof}[Proof of Proposition~\ref{prop:EProps}]
$E$ is lower semi-continuous as a rate function of a large deviations principle. Strict convexity and positivity away from $0$ are inherited from $e$. This shows \eqref{item:EProps1}. As for \eqref{item:EProps2}, by Jensen's inequality for any $0 \leq l < r \leq 1$ we have
\[
\int_l^r e(\gga^{\prime}(x)) d x 
	\geq 	(r-l) e \lb \inv{r-l} \int_l^r \gga^{\prime}(x) d x \rb
	=		(r-l) e \lb \frac{\gga(r) - \gga(l)}{r-l} \rb.
\]
This gives $E(\linear(\lbr (l, \gga(l)), (r, \gga(r)) \rbr)\1_{[l,r]}) \leq E(\gga \1_{[l,r]})$ and the proof is completed by summation.
\end{proof}

\begin{proof}[Proof of Proposition~\ref{Thm.Prop.Continuity}]
The first two parts are easy to verify. As for the third,
fix $\chi_0 < \infty$, $m_0 > 0$, $\gga \in \cL^0$, let $\gs, \wt{\gs} \in \gS_{\chi_0}(m_0)$
and set $d_0 = d(\gs, \wt{\gs})$.
We shall show that, once $d_0$ is small enough, we can explicitly construct
a $\wt{\gga} \in \cL^0$ satisfying \eqref{Item.MaxContinuity.1}, \eqref{Item.MaxContinuity.2}, \eqref{Item.MaxContinuity.3} with $\gep = \gep(\chi_0, m_0, d_0)$ independently of $\gga, \gs, \wt{\gs}$ and that
$\gep(\chi_0, m_0, d_0) \to 0$ as $d_0 \to 0$ for all $\chi_0$, $m_0$.
In what follows, we shall often omit the $\gs, \wt{\gs}$ subscript and instead  add tilde above quantities related to $\wt{\gs}$.
Let $\gi = (\gi_j)_{j=1}^{|\gi|} = I_{\gs}(\gga)$
and recall that the indices in $\gi$ are ordered according
to the $x$-coordinate of the indexed masses. This induces a piecewise decomposition of $\gga$:
\[
\gga = \sum_{j=1}^{|\gi|-1} \gga^j
	\quad ; \quad
\gga^j = \gga \1_{\lsp x^{\gi_j}, x^{\gi_{j+1}} \rb}.
\]
We will use this decomposition to construct $\wt{\gga}$. Formally 
for $j = 1, \dots, |\gi|-1$ set
\[
\wt{\gga}^j_{\pm}(x) =
    \wt{y}^{\gi_j} +
    a^j_{\pm} (\gga^j (x^{\gi_j} + b^j (x - \wt{x}^{\gi_j})) - y^{\gi_j}) +
    c^j_{\pm} (x - \wt{x}^{\gi_j}) \1_{\lsp \wt{x}^{\gi_j}, \wt{x}^{\gi_{j+1}} \rb},
\]
where
\[
\begin{array}{lcr}
a_{\pm}^j = \frac{\gD\wt{x}(\gi_j, \gi_{j+1}) \mp \gD\wt{y}(\gi_j, \gi_{j+1})}
                {\gD x(\gi_j, \gi_{j+1}) \mp \gD y(\gi_j, \gi_{j+1})} ;&
b^j = \frac{\gD x(\gi_j, \gi_{j+1})}{\gD \wt{x}(\gi_j, \gi_{j+1})} ;&
c_{\pm}^j = \pm \lb 1 - a_{\pm}^j b^j \rb
\end{array},
\]
and let $\wt{\gga} = \sum_{j=1}^{|\gi|-1} \wt{\gga}^j$ where
\be
\label{Eqn.DefOfGammaTildeJ}
    \wt{\gga}^j = \lbr \begin{array}{ll}
        \wt{\gga}^j_+   & \text{if} \quad \wt{a}(\gi_j, \gi_{j+1}) \geq a(\gi_j, \gi_{j+1}). \\
        \wt{\gga}^j_-   & \text{otherwise}.
    \end{array} \rpr
\ee
We now argue that $\wt{\gga}$ is in $\cL^0$ and satisfies \eqref{Item.MaxContinuity.1}, \eqref{Item.MaxContinuity.2}, 
\eqref{Item.MaxContinuity.3}.
Indeed, it is easy to verify that each piece $\wt{\gga}^j$ is supported on
$\lsp \wt{x}^{\gi_j}, \wt{x}^{\gi_{j+1}} \rb$ and satisfies
$\wt{\gga}^j \lb \wt{x}^{\gi_j} \rb = \wt{y}^{\gi_j}$ and
$\wt{\gga}^j \lb (\wt{x}^{\gi_{j+1}})^- \rb = \wt{y}^{\gi_{j+1}}$.
Also not difficult is
$a_{\pm}^j = 1 + o(1)$, $b^j = 1 + o(1)$, $c_{\pm}^j = o(1)$
as $d_0 \to 0$ and
\[
c^j_{\pm} = \frac{\wt{a}(\gi_j, \gi_{j+1}) - a(\gi_j, \gi_{j+1})}{1 \mp a(\gi_j, \gi_{j+1})},
\]
which implies that $\pm c_{\pm}^j \geq 0$ when it is used for $\wt{\gga}^j$ in \eqref{Eqn.DefOfGammaTildeJ}.
Then since
$\tfrac{d}{d x} \wt{\gga}^j_{\pm} = \tfrac{d}{d x} \gga^j + c^j_{\pm} \lb 1 \mp \tfrac{d}{d x} \gga^j \rb$
and $\labs \tfrac{d}{d x} \gga^j \rabs \leq 1$ on $\lb x^{\gi_j}, x^{\gi_{j+1}} \rb$, it follows that $\labs \frac{d}{d x} \wt{\gga}^j_{\pm} \rabs \leq 1$
on $\lb \wt{x}^{\gi_j}, \wt{x}^{\gi_{j+1}} \rb$ once $d_0$ is sufficiently small. This shows that $\wt{\gga} \in \cL^0$.

Now, since $\gga, \wt{\gga} \in \cL^0$, we have for all $j$ and
$x \in \lsp x^{\gi_j} \maxwith \wt{x}^{\gi_j}\ ,\: x^{\gi_{j+1}} \minwith \wt{x}^{\gi_{j+1}} \rb$
\begin{eqnarray*}
\labs \gga^j(x) - \wt{\gga}^j(x) \rabs
    & =     & \labs \wt{y}^{\gi_j} + a^j_{\pm} \lb \gga^j \lb x^{\gi_j} + b^j(x-\wt{x}^{\gi_j}) \rb - y^{\gi_j} \rb
                + c^j_{\pm} (x-\wt{x}^{\gi_j}) - \gga^j(x) \rabs    \\
    & =     & \labs \wt{y}^{\gi_j} + \gga^j(x) - y^{\gi_j} - \gga^j(x) + o(1) \rabs = o(1).
\end{eqnarray*}
as $d_0 \to 0$. Therefore ~\eqref{Item.MaxContinuity.1} follows from
\[
\normsup{\gga - \wt{\gga}}
    \leq  \big( \supone{1 \leq j \leq |\gi|-1} \normsup{ (\gga^j - \wt{\gga}^j)
                \1_{\lsp x^{\gi_j} \maxwith \wt{x}^{\gi_j}\ ,\ x^{\gi_{j+1}} \minwith \wt{x}^{\gi_{j+1}} \rb}} \big) +
                2 \supone{1 \leq j \leq |\gi|} |x^{\gi_j} - \wt{x}^{\gi_j}|   
	= o(1).
\]
For \eqref{Item.MaxContinuity.2} write
\[
\wt{\gs}(\wt{\gga}) \geq \sum_{j=1}^{|\gi|} \wt{v}^{\gi_j} \geq \sum_{j=1}^{|\gi|} (v^{\gi_j} + o(1))
    = \gs(\gga) + o(1).
\]
Finally \eqref{Item.MaxContinuity.3} comes from
\begin{eqnarray*}
E(\wt{\gga})
    & =     & \int_0^1 e \lb \frac{d}{dx} \wt{\gga} \rb d x =  \sum_{j=1}^{|\gi|-1} \int_{\wt{x}^{\gi_j}}^{\wt{x}^{\gi_{j+1}}} e \lb \frac{d}{d x} \wt{\gga} \rb d x   =     \sum_{j=1}^{|\gi|-1} \int_{\wt{x}^{\gi_j}}^{\wt{x}^{\gi_{j+1}}} e \lb \frac{d}{d x} \gga + o(1) \rb d x    \\
    & \leq  & \sum_{j=1}^{|\gi|-1} \frac{1}{b^j} \int_{x^{\gi_j}}^{x^{\gi_{j+1}}} e \lb \frac{d}{d x} \gga \rb d x + o(1)   =   E(\gga) + o(1),
\end{eqnarray*}
where we used the uniform continuity and boundedness of $e$ in $[-1, 1]$.
\end{proof}

\begin{proof}[Proof of Proposition~\ref{Thm.Prop.RWLD}]
We shall use a standard tilting argument.
Let $\mu$ be a probability measure under which $(\xi_i)_{i=1}^{\infty}$ are independent, symmetric $\pm 1$ random variables. Set $S_0 = 0$, $S_k = \sum_{i=1}^k \xi_i$ for $k=1, \dots$.
For any $\gl \in \bbR$, we denote by $\mu^{\gl}$ the exponential tilting of $\mu$, namely 
\[
	\mu^{\gl} (\xi_i = x) = \frac{e^{\gl x} \mu(\xi = x)}{\mu e^{\gl \xi_i}},
\]
where $\mu f$  denotes expectation of $f$ with respect to $\mu$.

It is easy to verify that
\begin{equation}
\label{Eqn.TiltdedProb}
	\mu^{\gl} (S_n = x) = \exp \{ \gl x - nL(\gl) \} \mu(S_n = x),
\end{equation}
where $L(\gl) = \log \lb \mu e^{\gl \xi_i} \rb$.
Moreover, all moments of $\xi_i$ under $\mu^{\gl}$ are $C^{\infty}$ as
a function of $\gl$ and in particular $\gl \mapsto \mu^{\gl} \xi_i = L^{\pr}(\gl)$
is increasing and tending to $\pm 1$ as $\gl \to \pm \infty$.

If $|x| < n$, we may choose $\gl = \gl(x)$ such that $\mu^{\gl} \xi_i =
L^{\pr}(\gl(x)) = \frac{x}{n}$,
in which case $e \lb \frac{x}{n} \rb = \gl(x) \frac{x}{n} - L(\gl(x))$
where $e(x) = \sup_{\gl \in \bbR} \{\gl x - L(\gl)\}$ and the latter is an implicit form of \eqref{eqn:LittleEDef} (see, for instance, Lemma 2.2.5 in ~\cite{DZ}). Then ~\eqref{Eqn.TiltdedProb} becomes
\begin{equation}
	\mu(S_n = x) = \exp \lbr -n e \lb \frac{x}{n} \rb \rbr \mu^{\gl(x)}(S_n = x).
\end{equation}
We may then use local Central Limit Theorem for $S_n$ under $\mu^{\gl}$,
which holds uniformly on any bounded set of $\gl$-s,
since in this case we have a uniform bound on moments of $\xi_i$.
It follows that for any $\gt < 1$,
there exists $C_1=C_1(\gt) > 0$, $C_2 = C_2(\gt) > 0$ such that
\[
	 C_1 \exp \lbr -n e \lb \frac{x}{n} \rb \rbr n^{-1/2}
		\leq
	\mu(S_n = x)
		\leq
	C_2 \exp \lbr -n e \lb \frac{x}{n} \rb \rbr n^{-1/2},
\]
for all $x$ such that $\mu(S_n = x) \neq 0$ and $\labs \frac{x}{n} \rabs < \gt$.
This implies
\begin{equation*}
-\inv{n} \log \mu(S_n = x) = e \lb \frac{x}{n} \rb + o(1).
\end{equation*}
Now fix $\chi_0 < \infty$, $m_0 > 0$ and let $\gs \in \gS_{\chi_0, n}(m_0)$ and $\gi \in \cI(\gs)$. Then
\begin{eqnarray*}
-\inv{n} \log \mu_n \lb \gi \subseteq I_{\gs}(s) \rb
	& 	=	& -\inv{n} \sum_{j=1}^{|\gi|-1} \log
				\mu \lb S_{n \gD x_{\gs} (\gi_j, \gi_{j+1})} = n \gD y_{\gs} (\gi_j, \gi_{j+1})	 \rb \\
	&	=	& \sum_{j=1}^{|\gi|-1} \gD x_{\gs}(\gi_j, \gi_{j+1}) \lb
				 e \lb a_{\gs}(\gi_j, \gi_{j+1}) \rb + o(1) \rb	\\
	&  	= 	& E \lb \gGa_{\gs} (\gi) \rb + o(1),
\end{eqnarray*}
uniformly as desired.

As for the second part, fix in addition $\gd > 0$ and let $\gs$, $\gi$ be as before. Then
\begin{eqnarray}
\nonumber
\lefteqn{\mu_n \lb \lpr \normsup{s-\gGa_{\gs}(\gi)} > \gd \rabs I_{\gs}(s) = \gi \rb}	 \\
\nonumber
	& \leq	& \sum_{j=1}^{|\gi|-1}
				\mu \big( \exists \, 0 \leq k \leq n \gD x_{\gs}(\gi_j, \gi_{j+1}) : \:
						|S_k - k\,a_{\gs}(\gi_j, \gi_{j+1})| > \gd n \ \big| 	 \\
\label{Eqn.RWLDPropEqn1}
	&		& \quad \quad 
					S_{n \gD x_{\gs}(\gi_j, \gi_{j+1})} = n \gD y_{\gs}(\gi_j, \gi_{j+1}), \,
					S_{n \gD x_{\gs}(\gi_j, i)} \neq n \gD y_{\gs}(\gi_j,i) ; \; \forall i \notin \gi
			  \big).
\end{eqnarray}
For the rest of the proof, we write $\gD x_j$,$\gD y_j$, $\gD a_j$
as a short for $\gD x_{\gs}(\gi_j, \gi_{j+1})$, 
$\gD y_{\gs}(\gi_j, \gi_{j+1})$, $\gD a_{\gs}(\gi_j, \gi_{j+1})$
and $\gD x_j(i)$,$\gD y_j(i)$ as a short for $\gD x_{\gs}(\gi_j, i)$,
$\gD y_{\gs}(\gi_j, i)$.
Choosing $\gl_j$ such that $\mu^{\gl_j} \xi_i = L^{\pr}(\gl_j) = a_j$ and
setting $\wt{S}_k = S_k - k\,a_j$, the $j$-th term in the above sum is equal to
\begin{eqnarray*}
\lefteqn{  \mu^{\gl_j} \lb \lpr \exists \ \gd n / 2 \leq k \leq n \gD x_j : \:
						| \wt{S}_k | > \gd n \ \rabs \  	
					\wt{S}_{n \gD x_j} = 0,  
					\wt{S}_{n \gD x_j(i)} \neq n \lb \gD y_j(i) -
						\gD x_j(i) a_j \rb
						; \; \forall i \notin \gi
			  \rb}	\\
	& \leq 	& \frac{\sum_{k=\gd n /2}^{n \gD x_j}
						\mu^{\gl_j} \lb | \wt{S}_k | > \gd n \rb}{
		\mu^{\gl_j} \lb \wt{S}_{n \gD x_j} = 0 \rb - 
					\sum_{i \notin \gi} \mu^{\gl_j} \lb
						\wt{S}_{n \gD x_j(i)} = n \lb \gD y_j(i) -
							\gD x_j(i) a_j \rb,
						\wt{S}_{n \gD x_j} = 0 \rb}.
\end{eqnarray*}
For the numerator, Cramer's Theorem implies
\[
	-\inv{n} \log \mu^{\gl_j} \lb | \wt{S}_k | > \gd n \rb \geq (\wt{e}^{\gl_j}(\gd)
		\minwith \wt{e}^{\gl_j}(-\gd)) + o(1),
\]
as $n \to \infty$, uniformly in the range of $k$, where $\wt{e}^{\gl}(\cdot)$ is Cramer's rate
function for $\xi_i - \mu^{\gl} \xi_i$ under $\mu^{\gl}$. It is easy to verify that $\wt{e}(x) \bydef \infone{\gl \in (-\infty, +\infty)} \wt{e}^{\gl}(x)$ is positive
away from $0$, whence there exists $\eta = \eta(\gd)$ such that
\begin{equation*}
\label{Eqn.RWLDPropEqn2}
-\inv{n} \log \big( \sum_{\gd n /2}^{n \gD x_j}
		\mu^{\gl_j} (| \wt{S}_k | > \gd n) \big) \geq
	\eta + o(1),
\end{equation*}
as $n \to \infty$ uniformly as desired.
On the other hand, by local CLT, we have
\begin{equation*}
\label{Eqn.RWLDPropEqn3}
	\mu^{\gl_j} \lb \wt{S}_{n \gD x_j} = 0 \rb = \gO(n^{-1/2})
\end{equation*}
and
\begin{equation*}
\label{Eqn.RWLDPropEqn4}
\mu^{\gl_j} \lb \wt{S}_{n \gD x_j(i)} = n \lb \gD y_j(i) -
	\gD x_j(i) a_j \rb,
	\wt{S}_{n \gD x_j} = 0 \rb	
	= O (n^{-1}),
\end{equation*}
as $n \to \infty$ uniformly in $\gs \in \gS_{\chi_0}(m_0)$, $\gi \in \cI(\gs)$ and $j$. Together, this implies that the denominator is $\gO(n^{-1/2})$ and the proof is complete.
\end{proof}

\subsection{Acknowledgements}
We would like to thank  G\'{e}rard Ben Arous for showing us the paper \cite{Ha07} and for fruitful discussions. For the latter we would also like to thank Chuck Newman. The research of both authors was supported in part by NSF Grant OISE-0730136. The research of the first author was also supported in part by NSF Grant DMS 0806180.

\end{document}